\documentclass[11pt]{article}
\usepackage{amsfonts}
\newcommand{\beq}{\begin{equation}}
\newcommand{\eeq}{\end{equation}}
\newtheorem{Theorem}{Theorem}[section]
\newtheorem{Lemma}[Theorem]{Lemma}
\newtheorem{Remark}[Theorem]{Remark}
\newtheorem{prop}{Proposition}[section]
\newcommand{\lbl}{\label}

\newcommand{\la}{\lambda}

\newcommand{\p}{\varphi}

\usepackage[english]{babel}
\usepackage[letterpaper,top=2cm,bottom=2cm,left=3cm,right=3cm,marginparwidth=1.75cm]{geometry}

\usepackage{amsmath}
\usepackage{amsfonts}
\usepackage[colorlinks=true, allcolors=blue]{hyperref}

\title{Multiplicity of Laplacian eigenvalues that can be represented by sum of two squares using number theory}
\author{Changfeng Zhou\\
Math major \\
Department of Mathematical Sciences\\
University of Cincinnati\\
Cincinnati, OH 45221\vspace{0.2in} \\
Taige Wang\\
Department of Mathematical Sciences\\
University of Cincinnati\\
Cincinnati, OH 45221}
\date{}

\begin{document}
\maketitle

\begin{abstract}
 In this article, we use results of Number Theory to prove the conjecture on eigenvalue problem of a 2D elliptic PDE proposed by P. Korman in his recent paper \cite{ref}: for any even integer $2k$, one can find an eigenvalue $N$ that can be represented as $N=a^{2}+b^{2}$, with integers $a\neq b$ and multiplicity $2k$, while for any odd integer $2k + 1$, one can find an integer $M$ that can be represented as $M=a^{2}+b^{2}$ with multiplicity $2k+1$. In addition, the manuscript gives the formula to find those $N$'s. 
\end{abstract}


\section{Introduction}

 \ \ \ \ In rectangular region $D=(0,\pi) \times (0,\pi)$ in $\mathbb{R}^2$, the following Laplacian's eigenvalue problem:

\begin{equation}\label{1-1}
	\begin{cases} \triangle u + \la u = 0 \quad \mbox{on $D$},&\\
	u = 0 \quad \mbox{on $\partial D$},&
	\end{cases}\end{equation}
 \noindent possesses only eigenvalues $\la _{nm}\equiv n^2+m^2$ with $n, m\in \mathbb{N}$, which corresponds to the eigenfunction set $\{\sin nx \sin my, \sin(mx)\sin(ny)\}$ (see e.g. L.C. Evans \cite{Evans} and Y. Pinchover et al \cite{PR}). These eigenvalues are not always be simple and can have multiple multiplicities. In \cite{ref}, P. Korman summarized a conjecture that the multiplicity of eigenvalues could be any positive integer. This statement was \emph{claimed but not proved} in his paper: for any positive even integer, one could find an eigenvalue $N$ as sum of two squares: $N=p^2+q^2, p\neq q\in \mathbb{N}$, so that $N$'s multiplicity is the given even integer; for any positive odd integer, one can also find eigenvalue $M= p^2+q^2$ meanwhile $M = r^2+r^2$ with different $p, q, r$ so that multiplicity of this eigenvalue is the given odd number. Such multiplicities differ due to eigenvalues can be represented by different forms of sums of two squares. In particular, $N = 2$ is simple as $2=1^2+1^2$, $N = 5$ is not simple but multiplicity is 2 as $5 = 1^2+2^2$, $N = 50$ is not simple but multiplicity is 3 as $50 = 1^2+7^2 = 5^2+5^2$. It is observed that the number of these presentations decides multiplicity great than 1 if the eigenvalue has more than one summation representations of $p^2+q^2$. It is worth to point out that this phenomenon of the conjecture is not true for 2D domains other than $(0, \pi)\times(0, \pi)$. For instance, if $D$ is a plate in $\mathbb{R}^2:\, x^2+y^2\le r^2$ , all eigenvalues have multiplicity two, except for the principal ones (which are simple; see paper \cite{ref} by P. Korman).\\

The aim of this paper is to use Number Theory to answer the aforementioned conjectured question:\\

\noindent{\bf Conjecture I.} \emph{Given any $n\in \mathbb{N}$, there exists an eigenvalue $N\in\mathbb{N}$ which has multiplicity of $n$.}\\

Or equivalently but more fundamentally: \\

\noindent{\bf Conjecture II.} \emph{Given $m\in \mathbb{N}$, there exists an integer $N\in\mathbb{N}$ that can be written as a sum of two squares in exactly $m$ different ways.} \\

In particular, we have 

\begin{eqnarray*}
 \begin{cases}
m=\frac{n}{2},\,\,{\rm if}\,\, n\,\,{\rm is\,\,even,}&\\\\
m=\frac{n-1}{2},\,\,{\rm if}\,\,n\,\,\rm{is\,\,odd.} &
\end{cases}
\end{eqnarray*}

We also want to point out that there have been classical Number Theory results on Conjecture II, the two squares problem (also referred to \emph{two squares theorem}, see, e.g., REU paper \cite{Bhaskar}). In present manuscript, we are dealing with such inverse problem as Conjecture I or II, and show there are even infinitely many ways to find such ${N}$.  \\

In the rest of paper, we will give out some observations and insights about the conjecture on the Gauss Plane $\mathbb{C}$ which lead to the main theorem in Section \ref{Sec2}; In Section \ref{Sec3}, we will present lemmas and proof of the theorem. 

\section{Preliminaries}\label{Sec2}
\setcounter{equation}{0}

If an integer N can be written as a sum of squares of two distinct integers $p, q>0$
\begin{equation}
    \label{eq1}
    N = p^2 + q^2=(p+iq)(p-iq),\quad  p\neq q,
\end{equation}
\noindent $N$ is the point with coordinates $(p, q)$ lying on the circle with radius $\sqrt{N}$ in the first quadrant of $\mathbb{C}$. However, there are three more points on same circle, for which (\ref{eq1}) is true: $(p, -q), (-p, q), (-p, -q)$. If we identify one of these points, the other points shall be known.\\

Let us show that there are exactly two ways to write $5^3 = 125$ in the form (\ref{eq1}). Since $5 = (1+2i)(1-2i), 5^3 = (1+2i)^3(1-2i)^3$ , a polynomial of degree $6$ with respect to imaginary unit $i$. We can factor $5^3 = (1+2i)^3\bar{z}$, where $\bar{z} = \overline{(1 + 2i)^3}$. Similarly, $5^3 = (1+2i)^2(1-2i)\bar{z}$, where we denote $\bar{z} = \overline{(1 + 2i)^2(1 - 2i)}$. Using this fashion ($\bar{z}$ denotes the complex conjugate of the first factor), the possibilities are

\begin{equation}
    \label{eq2}
    5^3 = (1+2i)^3\bar{z}= (1+2i)^2(1-2i)\bar{z}=(1+2i)(1-2i)^2\bar{z}=(1-2i)^3\bar{z}.
\end{equation}

Let $1+2i=\sqrt{5}e^{i\varphi}$, then $1-2i=\sqrt{5}e^{-i\varphi}$, and the polar angles of the first factors in (\ref{eq2}) are 

$$3\varphi,\,\,\varphi,\,\,-\varphi,\,\,-3\varphi.$$

We will not consider factors with negative $\varphi$, as they are complex conjugates of the factorizations already considered. This leaves us with

\begin{eqnarray*}
5^3&=&(1+2i)^3\bar{z}=(-11-2i)\bar{z},\\
5^3&=&(1+2i)^2(1-2i)\bar{z}=(5+10i)\bar{z},
\end{eqnarray*}

whence

\begin{eqnarray*}
5^3=11^2+2^2=10^2+5^2,    
\end{eqnarray*}
and there are no other possibilities to write $5^3$ as a sum of two squares.\\

Similarly, $5^5=(1+2i)^5(1-2i)^5$, a polynomial of degree $10$ in $i$. The different possibilities are exhausted by

$$5^5 = (1+2i)^5\bar{z}= (1+2i)^4(1-2i)\bar{z}=(1+2i)^3(1-2i)^2\bar{z},$$

and the polar angles of the first factors are

$$5\varphi,\,\,3\varphi,\,\,\varphi.$$

Thus

$$5^5=(1+2i)^5\bar{z}=(41-38i)\bar{z},$$
$$5^5=(1+2i)^4(1-2i)\bar{z}=(-55-10i)\bar{z},$$
$$5^5=(1+2i)^3(1-2i)\bar{z}=(25+50i)\bar{z}.$$

i.e., 

$$5^5=41^2+38^2=55^2+10^2=25^2+50^2,$$

and there are no other possibilities to write $5^5$ as a sum of two squares.\\

Moreover, consider $2\cdot 5^4=(1+i)(1-i)\cdot(1+2i)^4(1-2i)^4$. The different possibilities are exhausted by

$$2\cdot 5^4 = (1+i)(1+2i)^4\bar{z} = (1+i)(1+2i)^3(1-2i)\bar{z} = (1+i)(1+2i)^2(1-2i)^2\bar{z}$$
$$= (1-i)(1+2i)^4\bar{z} = (1-i)(1+2i)^3(1-2i)\bar{z} = (1-i)(1+2i)^2(1-2i)^2\bar{z}.$$

Since $1+i = \sqrt{2}e^{i\frac{\pi}{4}}$, $1-i = \sqrt{2}e^{-i\frac{\pi}{4}}$, the polar angles of the first factors are

$$\frac{\pi}{4}+4\p,\,\,\frac{\pi}{4}+2\p,\,\,\frac{\pi}{4},\,\,-\frac{\pi}{4}+4\p,\,\,-\frac{\pi}{4}+2\p,\,\,-\frac{\pi}{4}.$$

In these, complex numbers with angles 

$${\pi\over 4}+4\p,\,\,{\pi\over 4}+2\p,\,\, {\pi\over 4}$$

have real (imaginary) parts equal to imaginary (real) parts with ones with angles 

$$-{\pi\over 4}+4\p, -{\pi\over 4}+2\p, -{\pi\over 4}.$$ 

That is, for $z=Ae^{i(\frac{\pi}{4}+\p)}=a_1+b_1i, \tilde z=Ae^{i(-\frac{\pi}{4}+\phi)}=a_2+b_2i$, we have $|a_1| = |b_2|, |b_1| = |a_2|$, which leads to same presentation of two squares.\\

Thus 

\begin{eqnarray*}
2\cdot5^4 &=& (1+i)(1+2i)^4\bar{z} = (17-31i)\bar{z},\\
2\cdot5^4 &=& (1+i)(1+2i)^3(1-2i)\bar{z} = (-35+5i)\bar{z},\\
2\cdot5^4 &=& (1+i)(1+2i)^2(1-2i)^2\bar{z} = (25+25i)\bar{z} .
\end{eqnarray*}

Therefore,

$$2\cdot5^4 = 17^2+31^2 = 35^2+5^2 = 25^2+25^2.$$

Following this series of observations, we can claim:

\begin{prop}
\lbl{p2}
$(i)$ If $M=5^{2m-1}$, there are exactly $m$ different ways to represent $M=p^2+q^2,\ p\neq q$.\\

$(ii)$ If $M=2\cdot5^{2m}$ there are exactly $m$ different ways to represent $M=p^2+q^2,\ p\neq q$, and in addition, one way $M=r^2+r^2$.
\end{prop}
\textbf{Proof.}
$(i)$ Note $5^{2m-1}=(1+2i)^{2m-1}\cdot(1-2i)^{2m-1}$, a polynomial of degree $4m-2$ in $i$. The different possibilities are exhausted by

$$5^{2m-1}=(1+2i)^{2m-1}\bar{z}= (1+2i)^{2m-2}(1-2i)\bar{z}=\cdots=(1+2i)^m(1-2i)^{m-1}\bar{z},$$

and the polar angles of the first factors are

$$(2m-1)\varphi,(2m-3)\varphi,\cdots,3\varphi,\varphi.$$

We need to exclude the possibility that one of these angles is equal to $\pi$ minus another angle. If $3\varphi=\pi-\varphi$ (or $\varphi=\pi-3\varphi$), then $\varphi=\frac{\pi}{4}$, contradicting that tan$\: \varphi=2$. The possibility that one of these angles is equal to $\pi$ plus another angle is just a complex conjugate of the first possibility.\\

$(ii)$ Writing $2\cdot5^{2m}=(1+i)(1-i)\cdot(1+2i)^{2m}(1-2i)^{2m}$. Similar as previous example of $2\cdot 5^4$, different treating $1+i$\,\,(writing $1+i$ as the first factor or the second factor) gives identical representation of two squares, hence will not change the possibilities to writing $2\cdot5^{2m}$ as a sum of two squares compared to $5^{2m}$. And the different possibilities are exhausted by

$$2\cdot5^{2m}=(1+i)(1+2i)^{2m}\bar{z}=(1+i)(1+2i)^{2m-1}(1-2i)\bar{z}=\cdots=(1+i)(1+2i)^{m+1}(1-2i)^{m-1}\bar{z},$$

and 

$$2\cdot5^{2m}=(1+i)(1+2i)^m(1-2i)^m\bar{z}=(1+i)5^m\bar{z},$$

hence $r=5^m.$$\hspace*{\fill}$ $\Box$\\

Similarly as this proposition, one can consider other bases such $13, \cdots$, and their powers $13^{2m-1}, 2\times 13^{2m}, \cdots$, we might claim an abstract version of this proposition as main theorem, which is the conjecture in our manner, stated as follows:

\begin{Theorem}\label{th2-1}
Given the eigenvalue problem (\ref{1-1}), \\

(i) for  any even integer $n=2k, k\in\mathbb{N}$, there exists an eigenvalue $N\in\mathbb{N}$ whose multiplicity is $n$ or $2k$; in fact, 

$$N=p_1p_2^{k-1},\,\,p_1, p_2\equiv1({\rm mod}4).$$

(ii) for any odd integer $n=2k+1, k\in\mathbb{N}$, there exists an eigenvalue $M\in\mathbb{N}$ whose multiplicity is $2k+1$, and 

$$M=2p^{2k},\,\,p\equiv1({\rm mod4}).$$

Here $p_1, p_2, p$ are all prime numbers. 
\end{Theorem}

In the next section, we will prove the theorem by constructing the suitable $N$ form given the multiplicity $n$. 

\section{Proof of Theorem \ref{th2-1}}\label{Sec3}
\setcounter{equation}{0}

We would like to introduce two lemmas which are fundamental and well-known in number theory to establish our Theorem \ref{th2-1}. 

%\ Consider a circle with radius $\sqrt{N}$ in complex plan, the number of integer points in first quadrant on the circle equals to the number of ways to represent $N$ by sum of two squares($p^{2}+q^{2}$ and $q^{2}+p^{2}$ are considered different). In complex plane $N=p^{2}+q^{2}=(p+qi)(p-qi)$ where $p,q\in\mathbb{N}$. Let $z$ denote $p+qi$, then $N=z*\overline{z}$.

\begin{Lemma}\label{le3-1}
{\rm (}{\bf Fermat Theorem on sum of two squares} or {\bf Girard's Theorem}{\rm )} For an odd prime number $p$, it holds $p\equiv 1({\rm mod}4)$,\\

iff

$$p = x^2+y^2 = (x+yi)(x-yi), \,\,x,y\in \mathbb{N}. $$
\end{Lemma}

We can get this lemma by looking up number theory textbook \cite{Hardy}. Next, We would stress on two subsets of $\mathbb{N}$ and use them in the rest of paper. Denote 

\begin{eqnarray*}
\mathbb{P}=\{4k+1 | k\in \mathbb{N}\},\\
\mathbb{Q}=\{4k+3 | k\in \mathbb{N}\}.
\end{eqnarray*}

One can observe that all primes but $2$, lie in either $\mathbb{P}$ or $\mathbb{Q}$. 

\begin{Lemma}\label{le3-2}
(Prime Factorization) For any integer $N\in\mathbb{N}, N>1$, there exist powers $l, m_i, k_j\subset\mathbb{N}$ and $\{p_i\}\subset\mathbb{P}, \{q_j\}\subset\mathbb{Q}$ such that 

\begin{equation}\label{3-1}
N=2^l\prod_{i}p_i^{m_i}\prod_{j}q_j^{k_j}.
\end{equation}
\end{Lemma}

\noindent {\bf Proof.} It is an apparent reformulation from well-known prime factorization on $\mathbb{N}{\backslash}\{1\}$, due to the fact that all prime numbers are either 2, or in $\mathbb{P}$, or $\mathbb{Q}$. $\hspace*{\fill}$ $\Box$\\

%Factor $N$ into products of prime numbers $N=\prod p_i^{m_i}$. From Fermat's theorem on sums of two squares, for odd prime number $p$ 

%\begin{equation}
 %   p\equiv 1 (mod4)\Leftrightarrow p=x^{2}+y^{2}=(x+yi)(x-yi)
%\end{equation} 

Armed with these classical results from Number Theory, we would dig into the proof of Theorem \ref{th2-1}. We consider finding suitable $N$ which is equal to $a^2+b^2=(a+bi)(a-bi)$ on Gauss Plane $\mathbb{C}$, where $a, b\in\mathbb{N}$. We allow $N=c(a^2+b^2)$ and $c$ is a complete square integer. We articulate this structure of $N$ in terms of (\ref{3-1}) as the following key lemma:  

\begin{Lemma}\label{le3-3}
If $N=a^2+b^2$ with $a, b\in\mathbb{N}$, then $N$ takes form of (\ref{3-1}) with constraint $m_i\ge 1, k_j\ge 0$ and $l\ge 0$. Moreover, $k_j$ is even (including $k_j\equiv 0$) and $l$ is $1$ or even (including $l\equiv0$).  
\end{Lemma}

\noindent {\bf Proof.} We would split the construction of $N$ into three main steps:\\

\noindent {\bf Step 1: formulating $2^l$}\\

We need to divide discussion on whether a given integer $M$ whose multiplicity is $m$ and $2M$ are perfect square or not. If $M$ or $2M$ is perfect square, the multiplicity of $N = 2M$ will be $m + 1$ or $m - 1$; while neither of $M$ and $2M$ are a perfect square, $N = 2M$ has the same multiplicity as $M$. \\ 

{\bf Step 1.1: Case of $M$ where $M$ and $2M$ are not perfect square}\\

For every representation of $M = a^2 + b^2$ such that $a > b$, we have $N = 2M = (a+b)^2 + (a-b)^2$. And for every $N = 2M = c^2 + d^2$ such that $c > d$, we have $M = (\frac{c + d}{2})^2 + (\frac{c - d}{2})^2$. Note that $2 \mid N = c^2 + d^2$, so $\frac{c + d}{2}, \frac{c - d}{2} \in \mathbb{Z}$. Thus when neither $M$ and $2M$ are perfect square, $N = 2M$ has the same multiplicity as $M$.\\

{\bf Step 1.2: Case of $M$ where $M$ or $2M$ is a perfect square}\\

Consider an additional presentation of perfect square $M = a^2 + b^2$ where either $a = 0$ or $b = 0$. These particular cases could be neglected as one possible way to write $M$ as a sum of two squares. Without loss of generality, we set $b = 0$, as discussed in Step 1.1, there exists a corresponding presentation of $N = 2M = c^2 + d^2$ where $c = a+b$ and $d = a-b$ which will not be neglected as a presentation for $N$ since $c, d > 0$. At same time, all other presentations of $M$ have a corresponding representation of $N$. Consequently, the multiplicity of $N = 2M$ will be exactly $m + 1$ if $M$ has multiplicity $m$. \\

When $2M$ is a perfect square, we have $2^2 \mid 2M$, and $\frac{M}{2}$ is also a perfect square. Then there is a presentation of $M = r^2 + r^2$, i.e., $a = b(=r)$. And the corresponding representation of $2M$ will again be neglected since $d = a-b = 0$. Consequently, the multiplicity of $N = 2M$ will be exactly $m - 1$ if $M$ has multiplicity $m$. \\

\noindent {\bf Step 2: formulating $p_i^{m_i}$}\\

Since $p_i\in\mathbb{P}$, Lemma \ref{le3-1} leads to $p_i=z_i\overline{z_i}$, we can multiply $z_i$'s to form the complex factor $a+bi$, and the counterpart, $z_i$'s conjugates $\overline{z_i}$'s, form the second complex factor $a-bi$. Also, $\overline{z_i}$ can be in $a+bi$ while $z_i$ can be in the other factor.\\

Since each $p_i$ has such two ways to write: $z_i\overline{z_i}$ and $\overline{z_i}z_i$, $a+bi$ could have $0, 1, ...,m_i$ numbers of $z_i$'s, meanwhile the other factor $a-bi$ must absorb $\overline{z_i}$ as the fashion that $z_i$ has been in the first factor. We can tell that for every $p_i^{m_i}$,  $2p_i^{m_i}$ have $m_i +1$ possibilities with at least one $m_i\ge 1$. \\

In fact if we write each $p_i = z_i\overline{z_i} = \sqrt{p_i}e^{i\p_i}\times\sqrt{p_i}e^{-i\p_i}$,  

$$\prod_{i=1}^kp_i^{m_i} = \prod_{i=1}^k\sqrt{p_i}e^{i\sum_{i=1}^k(m_i-2t_i)\p_i}\cdot\overline{z}\qquad t_i \in \{0,1,\cdots,m_i\}\ for\ all\ 1\leq i\leq k$$

and the polar angles of the first factors are

$$\sum_{i=1}^km_i\p_i,\,\,(m_1-2)\p_1+\sum_{i=2}^km_i\p_i,\,\,(m_1-4)\p_1+\sum_{i=2}^km_i\p_i,\,\,\cdots,\,\,-m_1\p_1+\sum_{i=2}^km_i\p_i,\,\,\cdots$$

$$m_1\p_1+(m_2-2)\p_2+\sum_{i=3}^km_i\p_i,\,\,\cdots,\,\,-m_1\p_1-m_2\p_2+\sum_{i=3}^km_i\p_i,\,\,\cdots,\,\,\sum_{i=1}^k-m_i\p_i.$$

Following this fashion, if we have product $p_1^{m_1}p_2^{m_2}\cdots p_{k}^{m_k}$, then there is the number of polar angles: $(m_1+1)(m_2+1)\cdots(m_k+1)$ for $2p_1^{m_1}p_2^{m_2}\cdots p_{k}^{m_k}$. \\

Note if all $m_i$ are even integers, then $p_1^{m_1}p_2^{m_2}\cdots p_{k}^{m_k}$ is a perfect square, via above discussion, multiplicity of $p_1^{m_1}p_2^{m_2}\cdots p_{k}^{m_k}$ is 1 less than $2p_1^{m_1}p_2^{m_2}\cdots p_{k}^{m_k}$. And if $m_i$'s are odd integers for some $i$'s, $p_1^{m_1}p_2^{m_2}\cdots p_{k}^{m_k}$ has exactly same multiplicity as $2p_1^{m_1}p_2^{m_2}\cdots p_{k}^{m_k}$.\\

\noindent {\bf Step 3: formulating $q_j^{k_j}$ with even number $k_j$}\\

Owing to $q_j\in\mathbb{Q}$, it doesn't have product form with complex numbers. The power $k_j$ must be even, such that the result of $\prod_i p_i^{m_i}\times q_j^{k_j}$ can be represented as $(a+bi)(a-bi)$.\\

Due to this formulation, $q_j$'s have no impact on number of two-square representations, when multiplying with $\prod_i p_i^{m_i}$, and $k_j$ could be taken as $0$.$\hspace*{\fill}$ $\Box$\\

This time, we prove the main theorem getting a refined terse form of $N$ by using fundamentals provided by Lemma \ref{le3-3}. \\

\noindent{\bf Proof of Theorem 1.1.} Combining the steps we had in Lemma \ref{le3-3}, we successfully establish the existence of $N$ linked to sum of squared numbers in the theorem. The remaining of the proof is to clarify the steps reaching such $N$ which possesses the multiplicity of given $n$.\\

For $n=2k,\,\,k\in\mathbb{N}$, we adjust the $N$ in (\ref{3-1}) to be $N=p_1p_2^{k-1}$, where $p_1, p_2\in \mathbb{P}$. This $N$ has exactly $k$ different ways to represent $N$ as $N=a^2+b^2$, according to Lemma \ref{le3-3}.\\

For odd multiplicity $n=2k+1$, let $M=2p^{2k}$, where $p\in\mathbb{P}$, then there are exactly $k$ different ways to represent $M=a^2+b^2$ following fashion of Lemma \ref{le3-3}. In addition, $M$ can also be $2a^2$ for $a=p^{k}$. $\hspace*{\fill}$ $\Box$\\

We have the following remark since the infinite ways of selections of $p_1, p_2, p\in \mathbb{P}$:

\begin{Remark}
In the process in above discussion, the selection of $N$ is apparently not unique, and is infinitely many. 
\end{Remark}

 We would point that, the analogous conjecture in higher dimensions, is false. In fact, if dimensionality $d=3$, if $N$ has exactly one way that $N=p^2+q^2+r^2, p\neq q\neq r,$ multiplicity is 6; if $N = p^2 + q^2 + q^2,$ multiplicity is 3; else if $N=r^2+r^2+r^2$, multiplicity is 1. Therefore, the viable multiplicities are $3k$ or $3k+1, k\in \mathbb{N}$, hence we can not propose any nature integer $n$ being multiplicity unless $n=3k$ or $3k+1.$ We summarized this observation in the following proposition for $d\ge 3$: 

\begin{prop}
\lbl{p2}
For dimensinality $d$, eigenvalue $M=\sum_{i=1}^d a_i (a_i\neq a_j$ for different $i, j$) has multiplicity $n!$,  $M=\sum_{i=1}^d a_i\,\,(a_{i_1}=a_{i_2}=\cdots=a_{i_{d_i}}, d_1 + d_2 + \cdots d_m =d$, i.e., for some specific $i$'s, $a_i$'s are identical) has multiplicity $\frac{d!}{\prod_{i=1}^m(d_i)!}$, and $M=r^2+\cdots+r^2=nr^2$ is simple. Therefore, eigenvalues of elliptic problem (\ref{1-1}) in orthotope have only multiplicity of numbers lying in $\{dk,\,dk+1 | k\in\mathbb{N}\}$, a proper subset of $\mathbb{N}$.\\
\end{prop}

\section*{Acknowledgement}
Authors would fully thank for referees' precious suggestion and diligent work. \\

\noindent One of the authors, Taige Wang, would take this opportunity to thank Taft Research Center of University of Cincinnati by supporting him with Taft Awards. 

\bibliographystyle{plain}

\end{document}